\documentclass[aos,seceqn,nameyear,dvips]{arximspdf}
\usepackage{mathbh}
\usepackage{graphicx}

\doi{10.1214/10-AOS817}
\volume{38}
\issue{5}
\pubyear{2010}
\firstpage{3063}
\lastpage{3092}

\makeatletter

\newtheorem{theorem}{Theorem}[section]
\newtheorem{lemma}{Lemma}[section]
\newtheorem{corollary}{Corollary}[section]
\newtheorem{proposition}{Proposition}[section]
\newproclaim{example}{Example}
\newproclaim{remark}{Remark}

\newcommand{\R}{\mathbb{R}}

\newcommand{\PROB}{\mathbb{P}}

\newcommand{\Var}{\operatorname{Var}}
\newcommand{\A}{\mathcal{A}}
\newcommand{\B}{\mathcal{B}}
\newcommand{\C}{\mathcal{C}}
\newcommand{\scN}{\mathcal{N}}
\newcommand{\bX}{\mathbf{X}}
\newcommand{\bxi}{\bolds{\xi}}
\newcommand{\bzeta}{\bolds{\zeta}}

\newcommand{\diam}{\operatorname{diam}}

\makeatother

\begin{document}
\begin{frontmatter}

\title{On combinatorial testing problems\protect\thanksref{T1}}
\runtitle{Combinatorial testing problems}

\thankstext{T1}{Parts of this work were
done at the Bellairs Research Institute of McGill University and
the Saturna West Island Mathematics Center.}

\begin{aug}
\author[A]{\fnms{Louigi} \snm{Addario-Berry}\thanksref{t2}\ead[label=e1]{louigi@gmail.com}},
\author[B]{\fnms{Nicolas} \snm{Broutin}\ead[label=e2]{nicolas.broutin@inria.fr}},
\author[C]{\fnms{Luc} \snm{Devroye}\thanksref{t3}\ead[label=e3]{lucdevroye@gmail.com}}\\ and
\author[D]{\fnms{G\'abor} \snm{Lugosi}\corref{}\thanksref{t1}\ead[label=e4]{gabor.lugosi@gmail.com}}
\runauthor{Addario-Berry, Broutin, Devroye and Lugosi}
\affiliation{McGill University, INRIA, McGill University and
ICREA\break and Pompeu Fabra University}
\address[A]{L. Addario-Berry\\
Department of Mathematics and Statistics\\
McGill University\\
Montreal, H3A 2K6\\
Canada\\
\printead{e1}}
\address[B]{N. Broutin\\
INRIA Rocquencourt\\
78153 Le Chesnay\\
France \\
\printead{e2}\hspace*{11.53pt}}
\address[C]{L. Devroye\\
School of Computer Science\\
McGill University\\
Montreal, H3A 2A7\\
Canada\\
\printead{e3}}
\address[D]{G. Lugosi\\
ICREA\\
and Department of Economics\\
Pompeu Fabra University\\
Barcelona\\
Spain \\
\printead{e4}}
\end{aug}

\thankstext{t2}{Supported by an NSERC Discovery grant.}

\thankstext{t3}{Supported by NSERC.}

\thankstext{t1}{Supported by the Spanish Ministry of Science and
Technology Grant MTM2009-09063
and PASCAL2 Network of Excellence under EC Grant 216886.}

\received{\smonth{8} \syear{2009}}
\revised{\smonth{1} \syear{2010}}

%
\begin{abstract}
We study a class of hypothesis testing problems in which,
upon observing the realization of an $n$-dimensional Gaussian
vector, one has to decide whether the vector was drawn from a
standard normal distribution or, alternatively, whether there
is a subset of the components belonging to a certain given class
of sets whose elements have been ``contaminated,'' that is, have
a mean different from zero. We establish some general conditions
under which testing is possible and others under which testing
is hopeless with a small risk. The combinatorial and geometric
structure of the class of sets is shown to play a crucial role.
The bounds are illustrated on various examples.
\end{abstract}

%

\end{frontmatter}

\section{Introduction}\label{sec1}

In this paper, we study the following hypothesis testing problem
introduced by
\citet{ArCaHeZe08}. One observes an $n$-dimensional vector $\bX=
(X_1,\ldots,X_n)$. The null hypothesis $H_0$ is that the components of
$\bX$ are independent and identically distributed (i.i.d.) standard
normal random variables. We denote the probability measure and
expectation under $H_0$ by $\PROB_0$ and $\mathbb{E}_0$, respectively.

To describe the alternative hypothesis $H_1$, consider a class
$\C=\{S_1,\ldots,S_N\}$ of $N$ sets of indices such that $S_k
\subset\{1,\ldots,n\}$ for all $k=1,\ldots,N$. Under $H_1$, there
exists an $S \in\C$ such that
\[
X_i \mbox{ has distribution }
\cases{
\scN(0,1), & \quad if $i \notin S$, \cr
\scN(\mu,1), &\quad if $i \in S$,}
\]
where $\mu>0$ is a positive parameter. The components of $\bX$ are
independent under $H_1$ as well. The probability measure of $\bX$
defined this way by an $S \in\C$ is denoted by $\PROB_S$. Similarly, we
write $\mathbb{E}_S$ for the expectation with respect to $\PROB_S$.
Throughout, we will assume that every $S \in\C$ has the same
cardinality $|S|=K$.

A test is a binary-valued function $f\dvtx\R^n \to\{0,1\}$. If
$f(\bX)=0$ then we say that the test accepts the null hypothesis,
otherwise $H_0$ is rejected. One would like to design tests such that
$H_0$ is accepted with a large probability when $\bX$ is distributed
according to $\PROB_0$ and it is rejected when the distribution of
$\bX$ is $\PROB_S$ for some $S\in\C$. Following
\citet{ArCaHeZe08},
we consider the risk of a test $f$ measured by
%
%
\begin{equation}\label{eq:bayes_risk}
R(f) = \PROB_0\{f(\bX)=1\}
+ \frac{1}{N}\sum_{S\in\C} \PROB_S\{f(\bX)=0\} .
\end{equation}
This measure of risk corresponds to the view that, under the
alternative hypothesis, a set $S\subset\C$ is selected uniformly at
random and the components of $\bX$ belonging to $S$ have mean $\mu$. In
the sequel, we refer to the first and second terms on the right-hand
side of (\ref{eq:bayes_risk}) as the type I and type II errors,
respectively.

We are interested in determining, or at least estimating the value of
$\mu$ under which the risk can be made small. Our aim is to understand
the order of magnitude, when $n$ is large, as a function of $n$, $K$,
and the structure of $\C$, of the value of the smallest $\mu$ for which
risk can be made small. The value of $\mu$ for which the risk of the
best possible test equals $1/2$ is called \textit{critical}.

Typically, the $n$ components of $\bX$ represent weights over the $n$
edges of a given graph $G$ and each $S \in\C$ is a subgraph of $G$.
When $X_i \sim\scN(\mu,1)$ then the edge $i$ is ``contaminated'' and we
wish to test whether there is a subgraph in $\C$ that is entirely
contaminated.

In \citet{ArCaHeZe08}, two examples were studied in detail. In one
case, $\C$ contains all paths between two given vertices in a
two-dimensional grid and in the other $\C$ is the set of paths from
root to a leaf in a complete binary tree. In both cases, the order of
magnitude of the critical value of $\mu$ was determined.
\citet{ArCaDu09} investigate another class of examples in which
elements of $\C$ correspond to clusters in a regular grid. Both
\citet{ArCaHeZe08} and \citet{ArCaDu09} describe numerous
practical applications of problems of this type.

Some other interesting examples are when $\C$ is:
\begin{itemize}
\item
the set of all subsets $S\subset\{1,\ldots,n\}$ of size $K$;
\item
the set of all cliques of a given size in a complete graph;
\item
the set of all bicliques
(i.e., complete bipartite subgraphs)
of a given size in a complete bipartite graph;
\item
the set of all spanning trees of a complete graph;
\item
the set of all perfect matchings in a complete bipartite graph;
\item
the set of all sub-cubes of a given size of a binary hypercube.
\end{itemize}

The first of these examples, which lacks any combinatorial structure,
has been studied in the rich literature on multiple testing; see,
for example,
\citet{Ing99},
\citet{Bar02},
\citet{DoJi04} and the references therein.

As pointed out in \citet{ArCaHeZe08}, regardless of what $\C$ is, one may
determine explicitly the test $f^*$ minimizing the risk. It
follows from basic results of binary classification
that for a given vector $\mathbf{x}=(x_1,\ldots,x_n)$,
$f^*(\mathbf{x})=1$,
if and only if the ratio of the likelihoods of $\mathbf{x}$ under
$(1/N)\sum
_{S\in\C} \PROB_S$ and $\PROB_0$
exceeds $1$.
Writing
\[
\phi_0(\mathbf{x}) = (2\pi)^{-n/2} e^{-\sum_{i=1}^n x_i^2/2}
\]
and
\[
\phi_S(\mathbf{x}) = (2\pi)^{-n/2} e^{-\sum_{i\in S}(x_i-\mu
)^2/2-\sum
_{i\notin S} x_i^2/2}
\]
for the
probability densities of $\PROB_0$ and $\PROB_S$, respectively,
the likelihood ratio at $\mathbf{x}$ is
\[
L(\mathbf{x})
= \frac{{1/N}\sum_{S\in\C} \phi_S(\mathbf{x})}{\phi
_0(\mathbf{x})}
= \frac{1}{N} \sum_{S\in\C}
e^{\mu x_S - K\mu^2/2},
\]
where $x_S = \sum_{i\in S} x_i$. Thus, the optimal test is given by
\[
f^*(\mathbf{x}) = \mathbh{1}_{\{ L(\mathbf{x}) >1 \}} = \cases{
0, &\quad if $\displaystyle{\frac{1}{N} \sum_{S\in\C}
e^{\mu x_S - K\mu^2/2} \le1 }$, \vspace*{2pt}\cr
1, &\quad otherwise.}
\]
The risk of $f^*$ (often called the Bayes risk) may then be written
as
\begin{eqnarray*}
R^* & = & R^*_{\C}(\mu) = R(f^*) = 1 - \frac{1}{2} \mathbb{E}_0 |L(\bX
)-1| \\
& = & 1 - \frac{1}{2}
\int\biggl| \phi_0(\mathbf{x}) -\frac{1}{N}\sum_{S\in\C} \phi
_S(\mathbf{x}) \biggr|\, d\mathbf{x}.
\end{eqnarray*}
We are interested in the behavior of $R^*$ as a function of $\C$ and
$\mu$. Clearly, $R^*$ is a monotone decreasing function of $\mu$.
(This fact is intuitively clear and can be proved easily by
differentiating $R^*$ with respect to $\mu$.)
For $\mu$ sufficiently large, $R^*$ is close to zero
while for very small values of $\mu$, $R^*$ is near its maximum value~$1$,
indicating that testing is virtually impossible.
Our aim is to understand for what values of $\mu$ the transition occurs.
This depends on the combinatorial and geometric structure of the
class $\C$. We describe various general conditions in both directions
and illustrate them on examples.
\begin{remark*}[(An alternative risk measure)]
\citet{ArCaHeZe08}
also consider the risk measure
\[
\overline{R}(f) = \PROB_0\{f(\bX)=1\} + \max_{S\in\C} \PROB_S\{
f(\bX
)=0\}.
\]
Clearly, $\overline{R}(f) \ge R(f)$ and
when there is sufficient symmetry in $f$ and $\C$,
we have equality.
However, there are significant differences between the two
measures of risk. The alternative measure $\overline{R}$ obviously
satisfies the following monotonicity property: for a class
$\C$ and parameter $\mu>0$, let $\overline{R}{}^*_\C(\mu)$ denote
the smallest
achievable risk. If $\A\subset\C$ are two classes then for
any $\mu$, $\overline{R}{}^*_\A(\mu) \le\overline{R}{}^*_\C(\mu)$.
In contrast to this, the ``Bayesian'' risk measure $R(f)$
does not satisfy such a monotonicity property as is shown in
Section \ref{nonmonotone}.
In this paper, we focus on the risk measure $R(f)$.
\end{remark*}
%
%
\begin{remark*}
Throughout the paper we assume, for simplicity, that each set $S\in\C$
has the same cardinality $K$. We do this partly in order to avoid
technicalities that are not difficult but make the arguments less
transparent. At the same time, in many natural examples this condition
is satisfied. If $\C$ may contain sets of different size such that
all sets have approximately the same number of elements, then all arguments
go through without essential changes. However, if $\C$ contains sets
of very different size then the picture may change because large
sets become much easier to detect and small sets can basically
be ignored. Another approach to handle sets of different size, adopted
by 
\citet{ArCaDu09}, is to change the model
of the alternative hypothesis such that the level $\mu$ of contamination
is appropriately scaled depending on the size of the set $S$.
\end{remark*}

\subsection*{Plan of the paper}
The paper is organized as follows.
In Section \ref{simpletests}, we briefly discuss two suboptimal but
simple and general testing rules (the \textit{maximum test}
and the \textit{averaging test}) that imply sufficient conditions
for testability that turn out to be useful in many examples.

In Section \ref{lowerbounds}, a few general sufficient conditions
are derived for the impossibility of testing under symmetry assumptions
for the class.

In Section \ref{examples}, we work out several concrete examples,
including the class of all $K$-sets, the class of all cliques of
a certain size in a complete graph, the class of all perfect matchings
in the complete bipartite graph and the class of all spanning trees
in a complete graph.

In Section \ref{nonmonotone}, we show that, perhaps surprisingly,
the optimal risk is not monotone in the sense that
larger classes may be significantly easier to test than small ones,
though monotonicity holds under certain symmetry conditions.

In the last two sections of the paper, we use techniques developed
in the theory of Gaussian processes to establish upper and lower
bounds related to geometrical properties of the class $\C$.
In Section \ref{hellinger}, general lower bounds are derived in terms of
random subclasses and metric entropies of the class $\C$.
Finally, in Section \ref{typeone} we take a closer look at the type I
error of the optimal test and prove an upper bound that,
in certain situations, is significantly tighter than the
natural bound obtained for a general-purpose maximum test.

\section{Simple tests and upper bounds}
\label{simpletests}

As mentioned in the \hyperref[sec1]{Introduction}, the test $f^*$ minimizing the
risk is explicitly determined. However, the performance of this
test is not always easy to analyze. Moreover, efficient computation
of the optimal test is often a nontrivial problem though efficient algorithms
are available in many interesting cases. (We discuss computational
issues for the examples of Section \ref{examples}.)
Because of these reasons, it is often useful to consider simpler,
though suboptimal, tests. In this section, we briefly discuss two
simplistic tests, a test based on averaging and a test based on maxima.
These are often easier to analyze and
help understand the behavior of the optimal test as well. In many
cases, one of these tests turn out to have a near-optimal performance.

\subsubsection*{A simple test based on averaging}

Perhaps the simplest possible test is based on the fact that
the sum of the components of $\bX$ is zero-mean normal under $\PROB_0$
and has mean $\mu K$ under the alternative hypothesis. Thus, it
is natural to consider the \textit{averaging test}
\[
f(\mathbf{x}) = \mathbh{1}_{\{ \sum_{i=1}^n X_i > \mu K/2 \}}.
\]
\begin{proposition}
\label{average}
Let $\delta>0$. The risk of the averaging test $f$ satisfies
$R(f) \le\delta$
whenever
\[
\mu\ge\sqrt{\frac{8n}{K^2} \log\frac{2}{\delta}}.
\]
\end{proposition}
\begin{pf}
Observe that under $\PROB_0$, the statistic $\sum_{i=1}^n X_i$
has normal $\scN(0,n)$ distribution
while for each $S\in\C$, under $\PROB_S$,
it is distributed as $\scN(\mu K,n)$.
Thus, $R(f) \le2e^{-(\mu K)^2/(8n)}$.
\end{pf}

\subsubsection*{A test based on maxima}
\label{sec:max}

Another natural test is based on the fact that under the alternative
hypothesis for some $S\in\C$, $X_S= \sum_{i\in S} X_i$ is normal
$(\mu K, K)$.
Consider the \textit{maximum test}
\[
f(\mathbf{x}) = 1 \quad\mbox{if and only if}\quad
\max_{S\in\C} X_S\ge\frac{\mu K + \mathbb{E}_0 \max_{S\in\C} X_S}{2}.
\]
The test statistic $\max_{S\in\C} X_S$ is often referred to
as a \textit{scan statistic} and has been thoroughly studied
for a wide range of applications; see
\citet{GlNaWa01}. Here, we only need the following simple observation.
\begin{proposition}
\label{maxtest}
The risk of the maximum test $f$ satisfies
$R(f) \le\delta$ whenever
\[
\mu\ge\frac{\mathbb{E}_0 \max_{S\in\C} X_S}{K} + 2\sqrt{\frac{2}{K}
\log
\frac{2}{\delta}}.
\]
\end{proposition}

In the analysis, it is convenient to use the following
simple Gaussian concentration inequality; see
\citet{TsIbSu76}.
\begin{lemma}[(Tsirelson's inequality)]
\label{tsirelson}
Let $X=(X_1,\ldots,X_n)$ be an vector of $n$ independent
standard normal random variables.
Let $f:\R^n\to\R$ denote a Lipschitz function
with Lipschitz constant $L$ (with respect to the Euclidean distance).
Then for all $t >0$,
\[
\PROB\{ f(X) -\mathbb{E} f(X) \geq t \}
\leq
e^{- t^2/(2 L^2)} .
\]
\end{lemma}
\begin{pf*}{Proof of Proposition \ref{maxtest}} Simply
note that under the null hypothesis, for each $S\in\C$, $X_S$
is a zero-mean normally distributed random variable with variance $K=|S|$.
Since $\max_{S\in\C} X_S$ is a Lipschitz
function of $\bX$ with Lipschitz constant $\sqrt{K}$,
by Tsirelson's inequality, for all $t>0$,
\[
\PROB_0 \Bigl\{ \max_{S\in\C} X_S
\ge\mathbb{E}_0 \max_{S\in\C} X_S + t \Bigr\}
\le e^{-t^2/(2K)}.
\]
On the other hand, under $\PROB_S$ for a fixed $S\in\C$,
\[
\max_{S'\in\C} X_{S'}
\ge X_S \sim\scN(\mu K,K)
\]
and therefore
\[
\PROB_S \Bigl\{ \max_{S\in\C} X_S \le\mu K - t \Bigr\} \le e^{-t^2/(2K)},
\]
which completes the proof.
\end{pf*}

The maximum test is often easier to compute than the optimal test
$f^*$, though maximization is not always possible in polynomial time.
If the value of $\mathbb{E}_0 \max_{S\in\C} X_S$ is not exactly known, one
may replace it in the definition of $f$ by any upper bound
and then the same upper bound will appear in the performance bound.

Proposition \ref{maxtest} shows that
the maximum test is guaranteed to work whenever
$\mu$ is at least $\mathbb{E}_0 \max_{S\in\C} X_S/K + \mbox
{const.}/\sqrt{K}$.
Thus, in order to better understand the behavior of the maximum
test (and thus obtain sufficient conditions for the optimal test
to have a low risk), one needs to understand
the expected value of $\max_{S\in\C} X_S$ (under $\PROB_0$).
As the maximum of Gaussian processes have been studied
extensively, there are plenty of directly applicable results available
for expected maxima.
The textbook of
\citet{Tal05} is dedicated to this topic.
Here, we only recall some of the basic facts.

First, note that one always has
$\mathbb{E}_0 \max_{S\in\C} X_S \le\sqrt{2K \log N}$
but sharper bounds can be derived by chaining arguments; see
\citet{Tal05} for an elegant and advanced treatment.
The classical chaining bound of
\citet{Dud79} works as follows.
Introduce a metric on $\C$ by
\[
d(S,T) = \sqrt{\mathbb{E}_0 (X_S-X_T)^2}= \sqrt{d_H(S,T)},\qquad S,T
\in\C,
\]
where $d_H(S,T)= \sum_{i=1}^n \mathbh{1}_{\{ \mathbh{1}_{\{ i\in S \}
} \neq\mathbh{1}_{\{ i\in T \}} \}} $
denotes the Hamming distance. For $t>0$, let
$N(t)$ denote the $t$-covering number of $\C$ with respect to
the metric $d$, that is, the smallest number of open balls of radius
$t$ that cover $\C$. By Dudley's theorem, there exists a numerical constant
$C$ such that
\[
\mathbb{E}_0 \max_{S\in\C} X_S \le C \int_0^{\diam(\C)} \sqrt{ \log
N(t)} \,dt,
\]
where $\diam(\C)=\max_{S,T \in\C} d(S,T)$
denotes the diameter of the metric space $\C$.
Note that since
$|S| =K$ for all $S\in\C$, $\diam(\C) \le\sqrt{2K}$.
Dudley's theorem is not optimal but
it is relatively easy to use.
Dudley's theorem has been refined, based on ``majorizing measures,''
or ``generic chaining'' which gives sharp bounds; see, for example,
\citet{Tal05}.
%
\begin{remark*}[(The \textsc{vc} dimension)] \label{rmk:VC}
In certain cases, it is convenient to further bound Dudley's inequality
in terms of the \textsc{vc} dimension; see \citet{VaCh1971}. Recall that
the \textsc{vc} dimension
$V(\C)$ of $\C$ is the largest positive integer $m$ such that there
exists an $m$-element set $\{i_1,\ldots,i_m\} \subset\{1,\ldots,n\}$
such that for all $2^m$ subsets $A \subset\{i_1,\ldots,i_m\}$ there
exists an $S \in\C$ such that $S \cap\{i_1,\ldots,i_m\} =A$.
\citet{Hau95} proved that the covering numbers of $\C$ may be bounded
as
\[
N(t) \le e\cdot\bigl(V(\C)+1\bigr) \biggl(\frac{2en}{t^2}
\biggr)^{V(\C)},
\]
so by Dudley's bound,
\[
\mathbb{E}_0 \max_{S\in\C} X_S \le C \sqrt{V(\C) K \log n}.
\]
\end{remark*}
%
%
\begin{remark*}[(Tests based on symmetrization)]
An interesting alternative to the maximum test, proposed and
investigated by
\citet{DuRo06}
and
\citet{ArBlRo09a}, is
based on the idea that under the null hypothesis the distribution
of the vector $\bX$ does not change if the sign of each component
is changed randomly, while under the alternative hypothesis
the distribution changes. In \citet{DuRo06} and \citet{ArBlRo09a},
methods based on symmetrization and bootstrap are suggested and
analyzed. Such tests are meaningful and interesting in the setup
of the present paper as well and it would be interesting to
analyze their behavior.
\end{remark*}

\section{Lower bounds}
\label{lowerbounds}

In this section, we investigate conditions under which the risk of
any test is large. We start with a simple universal bound that
implies that regardless of what the class $\C$ is, small risk cannot
be achieved unless $\mu$ is substantially large compared to $K^{-1/2}$.

\subsubsection*{A universal lower bound}

An often convenient way of bounding the Bayes risk $R^*$
is in terms of the Bhattacharyya measure of affinity
[\citet{Bha46}]
\[
\rho=\rho_\C(\mu) = \tfrac{1}{2} \mathbb{E}_0 \sqrt{L(\bX)}.
\]
It is well known [see, e.g., \citet{DeGyLu95}, Theorem 3.1] that
\[
1 - \sqrt{1-4\rho^2} \le R^* \le2\rho.
\]
Thus, $2\rho$ essentially behaves as the Bayes error
in the sense that
$R^*$ is near $1$ when $2\rho$ is near $1$, and is small when $2\rho$
is small.
Observe that, by Jensen's inequality,
\[
2\rho=
\mathbb{E}_0 \sqrt{L(\bX)}
= \int\sqrt{\frac{1}{N} \sum_{S\in\C} \phi_S(\mathbf{x}) \phi
_0(\mathbf{x})} \,d\mathbf{x}
\ge\frac{1}{N} \sum_{S\in\C}
\int\sqrt{ \phi_S(\mathbf{x}) \phi_0(\mathbf{x})} \,d\mathbf{x}.
\]
Straightforward calculation shows that
for any $S\in\C$,
\[
\int\sqrt{ \phi_S(\mathbf{x}) \phi_0(\mathbf{x})} \,d\mathbf
{x}=e^{-\mu^2K/8}
\]
and therefore we have the following.
\begin{proposition}
\label{universal}
For all classes $\C$,
$R^* \ge1/2$ whenever
$\mu\le\sqrt{(4/K)}\times\break\sqrt{\log(4/3)}$.
\end{proposition}

This shows that no matter what the class $\C$ is, detection is
hopeless if $\mu$ is of the order of $K^{-1/2}$. This
classical fact goes back to
\citet{LeC70}.

\subsubsection*{A lower bound based on overlapping pairs}

The next lemma is due to
\citet{ArCaHeZe08}. For completeness, we recall their proof.
\begin{proposition}
\label{pairs}
Let $S$ and $S'$ be drawn independently, uniformly, at random from $\C$
and let $Z=|S\cap S'|$. Then
\[
R^* \ge1- \tfrac{1}{2} \sqrt{\mathbb{E} e^{\mu^2 Z} -1}.
\]
\end{proposition}
\begin{pf}
As noted above, by the Cauchy--Schwarz inequality,
\[
R^* = 1 - \tfrac{1}{2} \mathbb{E}_0 |L(\bX)-1|
\ge1 - \tfrac{1}{2} \sqrt{\mathbb{E}_0 |L(\bX)-1|^2}.
\]
Since $\mathbb{E}_0 L(\bX) = 1$,
\[
\mathbb{E}_0 |L(\bX)-1|^2=\Var_0(L(\bX)) = \mathbb{E}_0 [L(\bX)^2] -1.
\]
However, by definition $L(\bX)= \frac1 N \sum_{S\in\C} e^{\mu X_S -
K\mu^2/2}$, so
we have
\[
\mathbb{E}_0 [L(\bX)^2]
=
\frac{1}{N^2} \sum_{S,S'\in\C} e^{-K \mu^2} \mathbb{E}_0 e^{\mu(X_S+X_{S'})}.
\]
But
\begin{eqnarray*}
\mathbb{E}_0 e^{\mu(X_S+X_{S'})}
& = &
\mathbb{E}_0 [ e^{\mu\sum_{i\in S\setminus S'} X_i} e^{\mu\sum_{i\in
S'\setminus S} X_i}
e^{2\mu\sum_{i\in S\cap S'} X_i} ] \\
& = &
(\mathbb{E}_0 e^{\mu X} )^{2(K-|S\cap S'|)} (\mathbb{E}_0 e^{2\mu X}
)^{|S\cap S'|}
\\
& = &
e^{\mu^2 (K-|S\cap S'|)+2\mu^2|S\cap S'|},
\end{eqnarray*}
and the statement follows.
\end{pf}

The beauty of this proposition is that it reduces the problem to
studying a purely combinatorial quantity. By deriving upper
bounds for the moment generating function of the overlap $|S\cap S'|$ between
two elements of $\C$ drawn independently and uniformly at random,
one obtains lower bounds for the critical value of $\mu$.
This simple proposition turns out to be surprisingly powerful as
it will be illustrated in various applications below.

\subsubsection*{A lower bound for symmetric classes}

We begin by deriving some simple consequences of Proposition \ref{pairs}
under some general symmetry conditions on the class~$\C$. The following
proposition shows that the universal bound of
Proposition~\ref{universal} can be improved by a factor of
$\sqrt{\log(1+n/K)}$ for all sufficiently
symmetric classes.
\begin{proposition}
\label{symmetric}
Let $\delta\in(0,1)$.
Assume that $\C$ satisfies the following conditions of symmetry.
Let $S,S'$ be drawn independently and uniformly at random from $\C$.
Assume that:
\textup{(i)} the conditional distribution of $Z=|S\cap S'|$ given $S'$ is identical
for all values of $S'$;
\textup{(ii)} for any fixed $S_0\in\C$ and $i\in S_0$, $\PROB\{i\in S\}=K/n$.
Then $R^*\ge\delta$ for all $\mu$ with
\[
\mu\le\sqrt{\frac{1}{K}\log\biggl(1+\frac{4n(1-\delta)^2}{K} \biggr)}.
\]
\end{proposition}
\begin{pf}
We apply Proposition \ref{pairs}.
By the first symmetry assumption, it suffices to derive a suitable
upper bound for $\mathbb{E}[ e^{\mu^2 Z}]=\mathbb{E}[ e^{\mu^2 Z}| S']$
for an arbitrary $S'\in\C$.
After a possible relabeling, we may assume that $S'=\{1,\ldots,K\}$
so we can write $Z=\sum_{i=1}^K \mathbh{1}_{\{ i\in S \}}$.
By H\"older's inequality,
\begin{eqnarray*}
\mathbb{E}[ e^{\mu^2 Z}]
& = & \mathbb{E}\Biggl[ \prod_{i=1}^K e^{\mu^2 \mathbh{1}_{\{ i\in S \}}} \Biggr]
\\
& \le&
\prod_{i=1}^K (\mathbb{E}[ e^{K\mu^2 \mathbh{1}_{\{ i\in S \}}} ] )^{1/K}
\\ & = &
\mathbb{E}[ e^{K\mu^2 \mathbh{1}_{\{ 1\in S \}}} ]\qquad
\mbox{[by assumption (ii)]}
\\
& = & (e^{\mu^2K} -1 ) \frac{K}{n} +1.
\end{eqnarray*}
Proposition \ref{pairs} now implies the statement.
\end{pf}

Surprisingly, the lower bound of Proposition \ref{symmetric} is
close to optimal in many cases. This is true, in particular when
the class $\C$ is ``small,'' made precise in the following statement.
\begin{corollary}
\label{smallclass}
Assume that $\C$ is symmetric in the sense of Proposition \ref{symmetric}
and that it contains at most $n^\alpha$ elements where $\alpha>0$.
Then $R^*\ge1/2$ for all $\mu$ with
\[
\mu\le\sqrt{\frac{1}{K}\log\biggl(1+\frac{n}{K} \biggr)}
\]
and $R^*\le1/2$ for all $\mu$ with
\[
\mu\ge\sqrt{\frac{2\alpha}{K}\log n} + \sqrt{\frac{8 \log4}{K} }.
\]
\end{corollary}
\begin{pf}
The first statement follows from Proposition \ref{symmetric} while the
second from Proposition \ref{maxtest} and the fact that
$\mathbb{E}_0 \max_{S\in\C} X_S \le\sqrt{2K\log|\C|}$.
\end{pf}

The proposition above shows that for any small and sufficiently symmetric
class, the critical value of $\mu$ is of the order of
$\sqrt{(\log n)/K}$, at least if $K\le n^\beta$ for some $\beta\in(0,1)$.
Later, we will see examples of ``large'' classes for which
Proposition~\ref{symmetric} also gives a bound of the correct
order of magnitude.

\subsubsection*{Negative association}

The bound of Proposition \ref{symmetric}
may be improved significantly under an additional condition
of negative association that is satisfied in several interesting
examples (see Section \ref{examples} below). Recall that a collection
$Y_1,\ldots,Y_n$ of random variables is \textit{negatively associated}
if for any pair of disjoint sets $I,J\subset\{1,\ldots,n\}$ and
(coordinate-wise) nondecreasing functions $f$ and $g$,
\[
\mathbb{E}[ f(Y_i, i\in I) g(Y_j, j\in J) ]
\le\mathbb{E}[ f(Y_i, i\in I) ] \mathbb{E}[g(Y_j, j\in J) ].
\]
%
%
\begin{proposition}
\label{negass}
Let $\delta\in(0,1)$ and
assume that the class $\C$ satisfies the conditions
of Proposition \ref{symmetric}. Suppose that the labels are such
that $S'=\{1,2,\ldots,K\} \in\C$. Let $S$ be a randomly chosen
element of $\C$.
If the random variables
$\mathbh{1}_{\{ 1 \in S \}},\ldots,\mathbh{1}_{\{ K\in S \}}$ are
negatively associated,
then
$R^*\ge\delta$ for all $\mu$ with
\[
\mu\le\sqrt{\log\biggl(1+\frac{n\log(1+4(1-\delta)^2)}{K^2} \biggr)}.
\]
\end{proposition}
\begin{pf}
We proceed similarly to the proof of Proposition \ref{symmetric}. We have
\begin{eqnarray*}
\mathbb{E}[ e^{\mu^2 Z}]
& = & \mathbb{E}\Biggl[ \prod_{i=1}^K e^{\mu^2 \mathbh{1}_{\{ i\in S \}}} \Biggr]
\\ & \le&
\prod_{i=1}^K \mathbb{E}[ e^{\mu^2 \mathbh{1}_{\{ i\in S \}}} ]
\qquad\mbox{(by negative association)}
\\ & = &
\biggl( (e^{\mu^2} -1 ) \frac{K}{n} +1 \biggr)^K.
\end{eqnarray*}
Proposition \ref{pairs} and the upper bound above imply that $R^*$ at
least $\delta$ for all $\mu$ such that
\[
\mu\le\sqrt{\log\biggl(1+\frac{n ((1+4(1-\delta)^2)^{1/K}-1 )}{K} \biggr)}.
\]
The result follows by using $e^y\ge1+y$ with $y=K^{-1} \log
(1+4(1-\delta)^2)$. 
\end{pf}

\section{Examples}
\label{examples}

In this section, we consider various concrete examples and work out
upper and lower bounds for the critical range of $\mu$.

\subsection{Disjoint sets}
\label{canonical}

We start with the simplest possible case, that is,
when all $S\in\C$ are disjoint (and therefore $KN\le n$).
Fix $\delta\in(0,1)$.
Then, under $\PROB_0$, the $X_S$ are independent normal $(0,K)$ random
variables
and the bound $\mathbb{E}_0 \max_{S\in\C} X_S \le\sqrt{2K\log N}$
is close to being tight. By applying the maximum test $f$, we see that
$R^* \le R(f) \le\delta$
whenever
\[
\mu\ge\sqrt{\frac{2\log N}{K}} + 2\sqrt{\frac{2\log(2/\delta)}{K}}.
\]
To see that this bound gives the correct order of magnitude,
we may simply apply Proposition \ref{pairs}.
Here $Z$ may take two values:
\[
Z = K \qquad\mbox{with probability $1/N$}\quad\mbox{and}\quad
Z = 0 \qquad\mbox{with probability $1-1/N$.}
\]
Thus,
\[
\mathbb{E} e^{\mu^2Z} -1 = \frac{1}{N} (e^{\mu^2 K} -1 )
\le\frac{1}{N} e^{\mu^2 K}
\]
and therefore $R^* \ge\delta$ whenever
\[
\mu\le\sqrt{\frac{\log(4N(1-\delta)^2)}{K}}.
\]
So in this case the critical transition occurs when $\mu$ is of the order
of $\hspace*{-0.2pt}\sqrt{(1/K)\hspace*{-0.2pt}\log\hspace*{-0.25pt} N}$.
In Section \ref{hellinger}, we use this simple
lower bound to establish lower bounds for general classes $\C$ of sets.
Note that in this simple case one may directly analyze the risk
of the optimal test and obtain sharper bounds. In particular,
the leading constant in the lower bound is suboptimal.
However, in this paper our aim is to understand some general phenomena
and we focus on orders of magnitude rather than on nailing down sharp
constants.
%
%
\begin{remark*}[(Multiple hypothesis testing)]
Taking $S = \{\{1\},\ldots,\{n\}\}$,
$K=1$, and $N=n\ge2$ in the
above example,
we obtain a connection with multiple hypothesis testing.
In the latter, one tests ``$m_i = 0$'' against ``$m_i = \mu$''
for every $1 \le i \le n$,
and traditionally uses as test statistics $X_i , 1 \le i \le n$, to
build a
multiple testing procedure
(often rejecting all the hypotheses corresponding to large $X_i$),
see, for instance,
\citet{RoWo05},
\citet{ArBlRo09b}.
Such a procedure will reject the global null hypothesis
``$\forall i \in S, m_i= 0$'' if at least one of the
alternatives ``$m_i = \mu$'' is preferred.
The main difference with the approach taken in this paper
concerns the error rate. A multiple testing procedure
is generally calibrated to control measures of the type I error
like the family wise error rate or the
false discovery rate, while the tests defined in this paper
are designed to control the entire risk. Finally, in this example,
the subsets are disjoint which is the traditional framework
in multiple testing.
\end{remark*}

\subsection{$K$-sets}

Consider the example when $\C$ contains all sets $S \subset\{1,\ldots
,n\}$
of size $K$. Thus, $N={n\choose K}$.
As mentioned in the \hyperref[sec1]{Introduction}, this
problem is very well understood as sharp bounds and
sophisticated tests are available; see, for example,
\citet{Ing99},
\citet{Bar02},
\citet{DoJi04}. We include it for illustration
purposes only and we warn the reader that the obtained bounds
are not sharpest possible.

Let $\delta\in(0,1)$.
It is easy to see that the
assumptions of Proposition \ref{negass} are satisfied
[this follows, e.g., from Proposition 11 of
\citet{DuRa98}]
and therefore
$R^* \ge\delta$ for all
\[
\mu\le\sqrt{\log\biggl(1+\frac{n\log(1+4(1-\delta)^2)}{K^2} \biggr)}.
\]
This simple bound turns out to have the correct order of magnitude
both when $n \gg K^2$ [in which case it is of the order of
$\sqrt{\log(n/K^2 )}$] and when $n \ll K^2$ (when it is
of the order of $\sqrt{n/K^2}$).

This may be seen by considering the two simple tests described in
Section \ref{simpletests} in the two different regimes.
Since
\[
\frac{\mathbb{E}_0 \max_{S\in\C} X_S}{K} \le
\frac{\sqrt{2K\log{n\choose K}}}{K}
\le\sqrt{2\log\biggl(\frac{ne}{K} \biggr)},
\]
we see from Proposition \ref{maxtest} that when
$K=O (n^{(1-\varepsilon)/2} )$ for some fixed $\varepsilon>0$, then
the threshold value is of the order of $\sqrt{\log n}$.
On the other hand, when $K^2/n$ is bounded away from zero, then
the lower bound implied by Proposition \ref{negass} above is of the
order $\sqrt{n/K^2}$
and the averaging test provides a matching upper bound
by Proposition \ref{average}.

Note that in this example the maximum test is easy to compute
since it suffices to find the $K$ largest values among $X_1,\ldots,X_n$.

\subsection{Perfect matchings}

Let $\C$ be the set of all perfect matchings of the complete
bipartite graph $K_{m,m}$. Thus, we have $n=m^2$ edges and $N=m!$,
and $K=m$. By Proposition \ref{average} (i.e., the averaging test),
for $\delta\in(0,1)$, one has $R(f) \le\delta$
whenever $\mu\ge\sqrt{8\log(2/\delta)}$.

To show that this bound has
the right order of magnitude, we may apply Proposition \ref{negass}.
The symmetry assumptions hold obviously and the negative association
property follows from the fact that $Z=|S\cap S'|$ has the same
distribution as the number of fixed points in a random permutation.
The proposition implies that for all $m$,
$R^* \ge\delta$ whenever
\[
\mu\le\sqrt{\log\bigl(1+\log\bigl(1+4(1-\delta)^2\bigr)\bigr)}.
\]
Note that in this case the optimal test $f^*$ can be approximated
in a computationally efficient way. To this end, observe that
computing
\[
\frac{1}{N} \sum_{S\in\C} e^{\mu X_S}
= \frac{1}{m!} \sum_{\sigma} \prod_{j=1}^m e^{\mu X_{(j,\sigma(j))}}
\]
(where the summation is over all permutations of $\{1,\ldots,m\}$)
is equivalent to computing the permanent of an $m\times m$ matrix
with nonnegative elements. By a deep result of
\citet{JeSiVi04}, this may be done
by a polynomial-time randomized approximation.

\subsection{Stars}
\label{stars}

Consider a network of $m$ nodes in which each pair of nodes interacts.
One wishes to test if there is a corrupted node in the network
whose interactions slightly differ from the rest. This situation
may be modeled by considering the class of \textit{stars}.

\begin{figure}[b]

\includegraphics{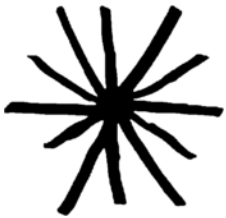}

\caption{A star [Vonnegut (\protect\citeyear{Von73})].}\label{asshole}
\end{figure}

A star is a subgraph of the complete graph $K_m$
which contains all $K=m-1$ edges containing a
fixed vertex (see Figure \ref{asshole}). Consider the
set $\C$ of all stars. In this setting, $n={m\choose2}$ and $N=m$.

%

In this case, we are in the situation of Corollary \ref{smallclass}
and Propositions \ref{symmetric} and~\ref{maxtest}
imply that
if $\C$ is the class of all stars in $K_m$ then for any $\varepsilon>0$,
\[
\lim_{m\to\infty} R^* = \cases{
0, &\quad if $\mu\ge\bigl(\sqrt{2}+\varepsilon\bigr)\sqrt{\dfrac{\log m}{m}}$,\vspace*{2pt}\cr
1, &\quad if $\mu\le(1-\varepsilon)\sqrt{\dfrac{\log m}{m}}$.}
\]

\subsection{Spanning trees}

Consider again a network of $m$ nodes in which each pair
of nodes interact. One may wish to test if there exists a
corrupted connected subgraph containing each node. This leads us to
considering the class of all spanning trees as follows.

Let $1,2,\ldots,n={m\choose2}$ represent the edges of
the complete graph $K_m$ and
let $\C$ be the set of all spanning trees of $K_m$.
Thus, we have $N=m^{m-2}$ spanning trees
and $K=m-1$. [See, e.g., \citet{Moon1970b}.] By Proposition \ref
{average}, the averaging
test has risk $R(f) \le\delta$
whenever $\mu\ge\sqrt{4\log(2/\delta)}$.

This bound is indeed of the right order. To see this, we may
start with Proposition \ref{pairs}. There are (at least) two
ways of proceeding. One is based on negative association.
Even though Proposition \ref{negass} is not applicable because
of the lack of symmetry in $\C$, negative association still holds.
In particular, by a result of
\citet{FeMi92}
[see also
\citet{GrWi04} and
\citet{BeLyPeSc01}],
if $S$ is a random uniform spanning tree of $K_m$, then
the indicators $\mathbh{1}_{\{ 1\in S \}},\ldots,\mathbh{1}_{\{ n\in
S \}}$ are
negatively associated. This means that, if $S$ and $S'$ are
independent uniform spanning trees and $Z=|S\cap S'|$,
\begin{eqnarray*}
\mathbb{E}[ e^{\mu^2 Z} ] & = &
\mathbb{E}\mathbb{E}\bigl[ e^{\mu^2 |S\cap S'|} \vert S' \bigr] \\
& = &
\mathbb{E}\mathbb{E}[ e^{\mu^2 \sum_{i\in S'}\mathbh{1}_{\{ i\in S \}}
} \vert
S' ] \\
& \le&
\mathbb{E}\prod_{i\in S'}\mathbb{E}[ e^{\mu^2 \mathbh{1}_{\{ i\in S \}} }
\vert S' ]\qquad
\mbox{(by negative association)} \\
& \le&
\mathbb{E}\prod_{i\in S'} \biggl( \frac{2}{m} e^{\mu^2 } +1 \biggr) \\
& = &
\biggl( \frac{2}{m} e^{\mu^2 } +1 \biggr)^{m-1} \\
& \le&
\exp(2e^{\mu^2 } ).
\end{eqnarray*}
This, together with Proposition \ref{pairs} shows that
for any $\delta\in(0,1)$,
$R^* \ge\delta$ whenever
\[
\mu\le\sqrt{\log\bigl(1+\tfrac1 2 \log\bigl(1+4(1-\delta)^2\bigr) \bigr)}.
\]
We note here that the same bound can be proved by a completely
different way that does not use negative association.
The key is to note that we may generate the two random spanning trees
based on $2(m-1)$ independent random variables $X_1,\ldots,X_{2(m-1)}$
taking values in
$\{1,\ldots,m-1\}$ as in
\citet{Ald90} [see also \citet{Broder1989}].
The key property we need is that if $Z_i$ denotes the number of common
edges in the two spanning trees
when $X_i$ is replaced by an independent copy $X_i'$ while keeping all
other $X_j$'s fixed, then
\[
\sum_{i=1}^{2(m-1)} (Z-Z_i)_+ \le Z
\]
(the details are omitted).
For random variables satisfying this last property,
an inequality of
\citet{BoLuMa00}
implies the sub-Poissonian bound
\[
\mathbb{E}\exp(\mu^2Z) \le\exp\bigl( \mathbb{E} Z(e^{\mu^2} -1) \bigr).
\]
Clearly, $\mathbb{E} Z = 2(n-1)/n \le2$, so essentially the same bound
as above
is obtained.

As the bounds above show, the computationally trivial average
test has a close-to-optimal performance. In spite of this, one may wish
to use the optimal test $f^*$. The ``partition function''
$(1/N) \sum_{S\in\C} e^{\mu X_S}$ may be computed by an algorithm
of
\citet{PrWi98}, who introduced a random sampling algorithm
that, given a graph with nonnegative weights $w_i$ over the edges,
samples a random spanning tree from a distribution such that the
probability of any spanning tree $S$ is proportional to
$\prod_{i\in S} w_i$.
The expected running time of the algorithm is bounded by the cover time
of an associated Markov chain that is defined as a random walk
over the graph in which the transition probabilities are proportional to
the edge weights. If $\mu$ is of the order of a constant (as in the critical
range) then the cover time is easily shown to be polynomial (with
high probability) as
all edge weights $w_i= e^{\mu^2 X_i}$ are roughly of the same order
both under the null and under the alternative hypotheses.

\subsection{Cliques}

Another natural application is the class of all cliques of a certain
size in a complete graph. More precisely, the random variables
$X_1,\ldots,X_n$
are associated with the edges of the complete graph
$K_m$ such that ${m\choose2}=n$ and let $\C$ contain all cliques
of size $k$. Thus, $K={k\choose2}$ and $N={m\choose k}$.
This case is more difficult than the class of $K$-sets discussed
above because negative association does not hold anymore.
(This may be easily seen by considering the indicator variables
of two adjacent edges both being in the randomly chosen clique.)
Also,
computationally the class of cliques is much more complex.
A related, well-studied model starts with the subgraph $K_m$ containing
each edge independently with probability $1/2$, as null hypothesis. The
alternative hypothesis is the same as the null hypothesis, except that
there is a clique of size $k$ on which each edge is independently
present with probability $p>1/2$. This is called the ``hidden clique''
problem (usually only the special case $p=1$ is considered). Despite
substantial interest in the hidden clique problem, polynomial time
detection algorithms are only known when $k=\Omega(\sqrt{n})$
[\citet{aks99}, \citet{fk00}]. We may obtain the hidden clique model from our
model by thresholding at weight zero (retaining only edges whose
normal random variable is positive), and so our model is easier for
testing than the hidden clique model. However, it seems likely that
designing an efficient test in the normal setting will be as difficult
as it has proved for hidden cliques.
It would be of interest to construct near-optimal tests that
are computable in polynomial time for larger values of $k$.

We have the following bounds for the performance of the optimal test.
It shows that when $k$ is a most of the order of $\sqrt{m}$,
the critical value of $\mu$ is of the order of $\sqrt{(1/k)\log(m/k)}$.
The proof below may be adjusted to handle larger values of $k$ as well
but we prefer to keep the calculations more transparent.
\begin{proposition}
\label{cliques}
Let $\C$ represent the class of all $N={m\choose k}$ cliques of
a complete graph $K_m$ and assume that $k \le\sqrt{m(\log2)/e}$.
Then:

\begin{longlist}
\item for all $\delta\in(0,1)$, $R^* \le\delta$ whenever
\[
\mu\ge
2\sqrt{\frac{1}{k-1} \log\biggl(\frac{me}{k} \biggr)}
+ 4\sqrt{\frac{\log(2/\delta)}{k(k-1)}};
\]
\item $R^* \ge1/2$ whenever
\[
\mu\le\sqrt{\frac{1}{k} \log\biggl(\frac{m}{2k} \biggr)}.
\]
\end{longlist}
\end{proposition}
\begin{pf}
(i) follows simply by a straightforward application of
Proposition \ref{maxtest} and the bound
$\mathbb{E}_0 \max_{S\in\C} X_S \le\sqrt{2 K \log N}$.

To prove the lower bound (ii), by Proposition \ref{pairs},
it suffices to show that if $S,S'$ are $k$-cliques drawn randomly and
independently
from $\C$ and $Z$ denotes the number of \textit{edges} in the intersection
of $S$ and $S'$, then $\mathbb{E}[ \exp(\mu^2Z) ] \le2$
for the indicated values of $\mu$.

Because of symmetry,
$\mathbb{E}[ \exp(\mu^2Z) ] = \mathbb{E}[ \exp(\mu^2Z) \vert S' ]$
for all $S'$ and therefore we might as well fix an arbitrary
clique $S'$. If $Y$ denotes the number of vertices in the clique
$S\cap S'$, then $Z={Y\choose2}$. Moreover, the distribution
of $Y$ is hypergeometrical with parameters $m$ and $k$.
If $B$ is a binomial random variable with parameters $k$ and $k/m$,
then since $\exp(\mu^2x^2/2)$ is a convex function of $x$, an
inequality of
\citet{Hoe63} implies that
\[
\mathbb{E}[ e^{\mu^2 Z} ]
= \mathbb{E}[ e^{\mu^2Y^2/2} ]
\le\mathbb{E}[ e^{\mu^2B^2/2} ].
\]
Thus, it remains to derive an appropriate upper bound for the moment
generating function
of the squared binomial. To this end, let $c>1$ be a parameter
whose value will be specified later. Using
\[
B^2 \le B \biggl(k \mathbh{1}_{\{ B> c{k^2}/{m} \}} + c\frac{k^2}{m} \biggr)
\]
and the Cauchy--Schwarz inequality, it suffices to show that
%
%
\begin{equation}
\label{product}
\mathbb{E}\biggl[\exp\biggl(\mu^2c\frac{k^2}{m}B \biggr) \biggr]
\cdot\mathbb{E}\bigl[\exp\bigl(\mu^2kB\mathbh{1}_{\{ B> c{k^2}/{m} \}} \bigr) \bigr]
\le4.
\end{equation}
We show that, if $\mu$ satisfies the condition of (ii),
for an appropriate choice of $c$,
both terms on the left-hand side are at most $2$.

The first term on the left-hand side of (\ref{product}) is
\[
\mathbb{E}\biggl[\exp\biggl(\mu^2c\frac{k^2}{m}B \biggr) \biggr]
= \biggl( 1+ \frac{k}{m} \biggl(\exp\biggl(\mu^2c\frac{k^2}{m} \biggr) -1 \biggr)
\biggr)^k,
\]
which is at most $2$ if and only if
\[
\frac{k}{m} \biggl(\exp\biggl(\mu^2c\frac{k^2}{m} \biggr) -1 \biggr)
\le2^{1/k} -1.
\]
Since $2^{1/k}-1 \ge(\log2)/k$, this is implied by
\[
\mu\le\sqrt{\frac{m}{ck^2} \log\biggl(1+ \frac{m\log2}{k^2} \biggr)}.
\]
To bound the second term on the left-hand side of (\ref{product}),
note that
\begin{eqnarray*}
\mathbb{E}\bigl[\exp\bigl(\mu^2kB\mathbh{1}_{\{ B> c{k^2/m} \}} \bigr) \bigr]
& \le&
1+ \mathbb{E}\bigl[\mathbh{1}_{\{ B> c{k^2/m} \}}\exp(\mu^2kB ) \bigr] \\
& \le&
1+ \biggl(\PROB\biggl\{B> c\frac{k^2}{m} \biggr\} \biggr)^{1/2}
(\mathbb{E}[\exp(\mu^2kB ) ] )^{1/2},
\end{eqnarray*}
by the Cauchy--Schwarz inequality,
so it suffices to show that
\[
\PROB\biggl\{B> c\frac{k^2}{m} \biggr\}\cdot
\mathbb{E}[\exp(\mu^2kB ) ]\le1.
\]
Denoting $h(x) = (1+x)\log(1+x)-x$, Chernoff's bound implies
\[
\PROB\biggl\{B> c\frac{k^2}{m} \biggr\}
\le\exp\biggl(-\frac{k^2}{m} h(c-1) \biggr).
\]
On the other hand,
\[
\mathbb{E}[\exp(\mu^2kB ) ]
= \biggl( 1 + \frac{k}{m}\exp(\mu^2k ) \biggr)^k,
\]
and therefore the second term on the left-hand side of (\ref{product})
is at most $2$ whenever
\[
1 + \frac{k}{m}\exp(\mu^2k )
\le\exp\biggl(\frac{k}{m} h(c-1) \biggr).
\]
Using $\exp(\frac{k}{m} h(c-1) ) \ge1+\frac{k}{m} h(c-1)$,
we obtain the sufficient condition
\[
\mu\le\sqrt{\frac{1}{k}\log h(c-1)}.
\]
Summarizing, we have shown that $R^* \ge1/2$ for all $\mu$
satisfying
\[
\mu\le
2\cdot\min\Biggl(\sqrt{\frac{1}{k}\log h(c-1)},
\sqrt{\frac{m}{ck^2} \log\biggl(1+ \frac{m\log2}{k^2} \biggr)} \Biggr).
\]
Choosing
\[
c= \frac{m}{k}\frac{\log(m/k)}{\log(m\log2/k^2)}
\]
[which is greater than $1$ for $k\le\sqrt{m(\log2)/e}$],
the second term on the right-hand side is at most
$\sqrt{(1/k)\log(m/k)}$. Now observe that since
$h(c-1)=c\log c -c+1$ is convex, for any $a>0$,
$h(c-1) \ge c\log a -a +1$. Choosing
$a=\frac{\log(m/k)}{\log(m\log2/k^2)}$,
the first term is at least
\[
\sqrt{\frac{1}{k}
\log\biggl(\frac{m}{k}-\frac{\log(m/k)}{\log(m\log2/k^2)} \biggr)}
\ge
\sqrt{\frac{1}{k} \log\biggl(\frac{m}{2k} \biggr)},
\]
where we used the condition that $m\log2/k^2 \ge e$ and
that $x \ge2\log x$ for all $x>0$.
\end{pf}
\begin{remark*}[(A related problem)]
A closely related problem arising in
the exploratory analysis of microarray data
[see
\citet{ShWePeNo09}]
is when each member of $\C$ represents the $K$ edges of a
$\sqrt{K}\times\sqrt{K}$ biclique
of the complete bipartite graph $K_{m,m}$ where $m=\sqrt{n}$.
(A biclique is a complete bipartite subgraph of $K_{m,m}$.)
The analysis and the bounds are completely analogous to the
one worked out above, the details are omitted.
\end{remark*}

\section{On the monotonicity of the risk}
\label{nonmonotone}

Intuitively, one would expect that the testing problem becomes harder
as the class $\C$ gets larger. More precisely, one may expect
that if $\A\subset\C$ are two classes of subsets of $\{1,\ldots,n\}$,
then $R^*_\A(\mu)\le R^*_\C(\mu)$ holds for all $\mu$.
The purpose of this section is to
show that this intuition is wrong in quite a strong sense
as not only such general monotonicity property does not hold for the risk,
but there are classes $\A\subset\C$ for which $R^*_\A(\mu)$ is arbitrary
close to~$1$ and $R^*_\C(\mu)$ is arbitrary close to $0$ for the
same value of $\mu$.

However, monotonicity does hold if the class $\C$ is
sufficiently symmetric.
Call a class $\C$ \textit{symmetric} if for the optimal test
\[
f^*_\C(\mathbf{x}) = \mathbh{1}_{\{ (1/N) \sum_{S\in\C} \exp(\mu
\sum_{i\in S} x_i) \ge\exp(K\mu^2/2) \}},
\]
the value of $\PROB_T\{f^*_\C(\bX)=0\}$ is the same
for all $T\in\C$.
Note that several of the examples discussed in Section \ref{examples}
satisfy the symmetry assumption, such as the classes of $K$-sets,
stars, perfect matchings, and cliques. However, the class of spanning
trees is not symmetric in the required sense.
\begin{theorem}
\label{symmetricclass}
Let $\C$ be a symmetric class of subsets of $\{1,\ldots,n\}$.
If $\A$ is an arbitrary subclass of $\C$,
then for all $\mu>0$,
$R^*_\A(\mu) \le R^*_\C(\mu)$.
\end{theorem}
\begin{pf}
In this proof, we fix the value of $\mu>0$
and suppress it in the notation.
Recall the definition of the alternative risk measure
\[
\overline{R}_\C(f) = \PROB_0\{f(\bX)=1\} + \max_{S\in\C} \PROB
_S\{f(\bX)=0\},
\]
which is to be contrasted with our main risk measure
\[
R_\C(f) = \PROB_0\{f(\bX)=1\} + \frac{1}{N}\sum_{S\in\C} \PROB
_S\{f(\bX
)=0\}.
\]
The risk $\overline{R}$ is obviously monotone in the sense that if $\A
\subset
\C$
then for every $f$, $\overline{R}_\A(f) \le\overline{R}_\C(f)$.
Let $\overline{f}{}^*_\C$ and $f^*_\C$ denote the optimal tests with
respect to both measures of risk.

First, observe that if $\C$ is symmetric, then
$\overline{R}_\C(f^*_\C) = R_\C(f^*_\C)$. But since $R_\C(f) \le
\overline
{R}_\C
(f)$ for every $f$,
we have
\[
\overline{R}_\C(\overline{f}{}^*_\C)\le\overline{R}_\C(f^*_\C)
= R_\C(f^*_\C) \le R _\C(\overline{f}{}^*_\C) \le\overline{R}_\C
(\overline{f}{}^*_\C).
\]
This means that all inequalities are equalities and, in particular,
$\overline{f}{}^*_\C= f^*_\C$.

Now if $\A$ is an arbitrary subclass of $\C$, then
\[
R^*_\C= R_\C(f^*_\C) = \overline{R}_\C(\overline{f}{}^*_\C)
\ge\overline{R}_\A(\overline{f}{}^*_\C)
\ge R_\A(\overline{f}{}^*_\C) \ge R_\A(f^*_\A) = R^*_\A,
\]
which completes the proof.
\end{pf}
\begin{theorem}
For every $\varepsilon\in(0,1)$ there exist $n$, $\mu$,
and classes $\A\subset\C\subset\{1,\ldots,n\}$
such that $R^*_\A(\mu) \ge1- \varepsilon$ and $R^*_\C(\mu) \le
2\varepsilon$.
\end{theorem}
\begin{pf} We work with $L_1$ distances.
For any class $\mathcal L$, denote $\phi_{\mathcal L}(\mathbf{x}) =
\frac
{1}{N}\times\sum_{S\in\mathcal L} \phi_S(\mathbf{x})$.
Recall that
\[
R^*_{\mathcal L}(\mu)
= 1 - \frac{1}{2}
\int| \phi_0(\mathbf{x}) -\phi_{\mathcal L}(\mathbf{x}) | \,d\mathbf{x}.
\]

Given $\varepsilon$, we fix an integer $K=K(\varepsilon)$ large enough that
$K+1 \geq1/\varepsilon$ and that
\[
\sqrt{\frac{\log(4(K+1)\varepsilon^2)}{K+1}}
\geq\sqrt{\frac{8}{K} \log\biggl(\frac{2}{\varepsilon} \biggr)},
\]
and let $n=n(\varepsilon)=(K+1)^2$.
We let $\mathcal A$ consist of $K+1$ disjoint subsets of $\{1,\ldots
,n\}
$, each of size $K+1$.
We let $\mathcal B$ consist of all sets of the form $\{1,\ldots,K,i\}$,
where $i$ ranges from $K+1$ to $n$,
and assume $\mathcal A$ has been chosen so that $\mathcal A \cap
\mathcal B= \varnothing$. We then let
$\mathcal C = \mathcal A \cup\mathcal B$. We take
%
\[
\mu= \sqrt{\frac{\log(4(K+1)\varepsilon^2)}{K+1}},
\]
so that, as seen in Section \ref{canonical}, we have $R^*_{\mathcal
A}(\mu) \geq1-\varepsilon$.
We will require an upper bound on $R^*_{\mathcal B}(\mu)$,
which we obtain by considering the averaging test on variables
$1,\ldots,K$,
\[
f(\mathbf{x}) = \mathbh{1}_{\{ \sum_{i=1}^K x_i \geq({\mu K})/{2} \}}.
\]
Just as in Proposition \ref{average}, we have $R(f) \leq\varepsilon$ whenever
$\mu\geq\sqrt{\frac{8}{K} \log(\frac{2}{\varepsilon} )}$, which is
indeed the case
by our choices of $\mu$ and $K$. It follows that $R^*_{\mathcal B}(\mu)
\leq\varepsilon$.
We remark that
\[
\int| \phi- \phi_{\A} | = 2-2R^*_{\mathcal A}(\mu) \leq2\varepsilon.
\]
%
We let $M=|\mathcal B|=(K+1)^2-K$; then $N=|\mathcal
C|=M+K+1=(K+1)^2+1$, and note
\begin{eqnarray*}
\int| \phi- \phi_{\C} |
&= & \int\biggl| \phi- { (K+1) \phi_{\A} + M \phi_{\B} \over N} \biggr| \\
&= & \int\biggl| { (K+1) (\phi- \phi_{\A} ) + M (\phi- \phi_{\B} )
\over
N } \biggr| \\
&\ge& {M \over N} \int| \phi- \phi_{\B} |
- {(K+1) \over N} \int| \phi- \phi_{\A} | \\
&\ge& (1- \varepsilon) \int| \phi- \phi_{\B} | - 2\varepsilon^2 \\
& = & (1-\varepsilon)\bigl(2-2R^*_{\mathcal B}(\mu)\bigr) - 2\varepsilon^2 \\
& \ge& 2-4\varepsilon.
\end{eqnarray*}
Thus, $R^*_{\mathcal C}(\mu) \leq2\varepsilon$.
\end{pf}

Observe that nonmonotonicity of the Bhattacharyya affinity also follows
from the same argument. To this end, we may express
$\rho_\C(\mu) =\break \frac12 \int\sqrt{\phi_0(\mathbf{x}) \phi
_S(\mathbf{x})} \,d\mathbf{x}$
in function of the Hellinger distance
\[
H(\phi_0,\phi_\C) = \sqrt{ \int\bigl( \sqrt{\phi_0(\mathbf{x})} -
\sqrt{\phi_\C(\mathbf{x}
)} \bigr)^2 \,d\mathbf{x}}
\]
as $\rho_\C(\mathbf{x}) = \frac1 2- \frac1 4 H(\phi_0,\phi_\C)^2$.
Recalling [see, e.g.,
\citet{DeGy85}, page 225] that
\[
H(\phi_0,\phi_\C)^2 \le\int| \phi_0(\mathbf{x}) -\phi_\C
(\mathbf{x}) | \le2H(\phi
_0,\phi_\C),
\]
we see that the same example as in the proof above, for $n$ large enough,
shows the nonmonotonicity of the Bhattacharyya affinity as well.

\section{Lower bounds on based random subclasses and metric entropy}
\label{hellinger}

In this section, we derive lower bounds for the Bayes risk
$R^* = R^*_\C(\mu)$.
The bounds are in terms of some geometric features of the class $\C$.
Again, we treat $\C$ as a metric space equipped with
the canonical distance
$d(S,T) = \sqrt{\mathbb{E}_0 (X_S-X_T)^2}$ [i.e., the square root of the
Hamming distance $d_H(S,T)$].

For an integer $M \le N$, we define a real-valued parameter
$t_\C(M)>0$ of the class $\C$ as follows. Let $\A\subset\C$
be obtained by choosing $M$ elements of $\C$ at random, without
replacement.
Let the random variable $\tau$ denote the smallest distance
between elements of $\A$ and let $t_\C(M)$ be a median
of $\tau$.
\begin{theorem}
\label{randomclassbound}
Let $M\le N$ be an integer. Then for any class $\C$,
\[
R_\C^* \ge1/4,
\]
whenever
\[
\mu\le\min\Biggl(\sqrt{\frac{\log(M/16)}{K}}, \frac{8\log(\sqrt
{3}/8)}{\sqrt{K-t_\C(M)^2/2}} \Biggr).
\]
\end{theorem}

To interpret the statement of the theorem, note that
\[
K-\tau^2/2 = {\mathop{\max_{S,T\in\A}}_{S\neq T}}
|S\cap T|
\]
is the largest overlap between any pair of elements of $\A$.
Thus, just like in Proposition \ref{pairs}, the distribution of
the overlap between random elements of $\C$ plays a key role
in establishing lower bounds for the optimal risk. However,
while in Proposition \ref{pairs} the moment generating function
$\mathbb{E}\exp(\mu^2|S\cap T|)$ of
the overlap between two random elements
determines an upper bound for the critical value of $\mu$,
here it is the median of the largest overlap between many
random elements that counts. The latter seems to carry more information
about the fine geometry of the class. In fact, invoking a simple
union bound, upper bounds for $\mathbb{E}\exp(\mu^2|S\cap T|)$
may be used together with Theorem \ref{randomclassbound}.

In applications, often it suffices to consider the following special
case.
\begin{corollary}
\label{medianzero}
Let $M\le N$ be the largest integer for which zero is a median of
${\max_{S,T\in\A,S\neq T}} |S\cap T|$
where $\A$ is a random subset of $\C$ of size $M$
[i.e., $t_\C(M)^2 = 2K$]. Then $R^*_\C(\mu) \ge1/4$
for all $\mu\le\sqrt{\log(M/16)/K}$.
\end{corollary}
\begin{example*}[(Sub-squares of a grid)]
To illustrate the corollary, consider the following example which is
the simplest in a family of problems investigated by
\citet{ArCaDu09}: assume that
$n$ and $K$ are both perfect squares and that the indices
$\{1,\ldots,n\}$ are arranged in a $\sqrt{n} \times\sqrt{n}$ grid.
The class $\C$ contains all $\sqrt{K} \times\sqrt{K}$ sub-squares.
Now if $S$ and $T$ are randomly chosen elements of $\C$ (with or
without replacement) then, if $(K+1)^2\le2\sqrt n$,
\[
\PROB\{ |S\cap T| \neq0\} \ge\frac{(\sqrt n - 2K)^2}{(\sqrt n-K+1)^2}
\cdot\frac{K}{(\sqrt n -K +1)^2} \ge\frac{K}{n}
\]
and therefore
\[
\PROB\Bigl\{ {\mathop{\max_{S,T\in\A}}_{S\neq T}} |S\cap
T| =
0 \Bigr\}
= 1 - \PROB\Bigl\{ {\mathop{\max_{S,T\in\A}}_{S\neq T}}
|S\cap T| > 0 \Bigr\}
\ge1 - M^2 \frac{K}{n},
\]
which is at least $1/2$ if $M \le\sqrt{n/(2K)}$ in which case
$t_\C(M)^2 = 2K$. Thus, by Corollary \ref{medianzero},
$R^*_\C(\mu) \ge1/4$
for all $\mu\le\sqrt{\log(n/(512 K))/(2K)}$. This bound is of the
optimal order of magnitude as it is easily seen by an application
of Proposition~\ref{maxtest}.
\end{example*}


In some other applications, a better bound is obtained if
some overlap is allowed. A case in point is the example of stars
from Section \ref{stars}. In that case, any two elements of $\C$
overlap but
by taking $M=N(=m)$, we have
$K-t_\C(M)^2/2=1$, so Theorem \ref{randomclassbound} still implies
$R^*_\C(\mu) \ge1/4$ whenever
$\mu\le\sqrt{(1/K)\log(m/16)}$.

The main tool of the proof of Theorem \ref{randomclassbound}
is Slepian's lemma which we recall here [\citet{Slepian1962}].
[For this
version, see
\citet{LeTa91}, Theorem~3.11.]
\begin{lemma}[(Slepian's lemma)]
\label{slepian}
Let $\bxi=(\xi_1,\ldots,\xi_N),\bzeta=(\zeta_1,\ldots,\zeta_N)
\in\R^N$
be zero-mean
Gaussian vectors such that for each $i,j=1,\ldots,N$,
\[
\mathbb{E}\xi_i^2 = \mathbb{E}\zeta_i^2 \qquad\mbox{for each
$i=1,\ldots,N$}\quad\mbox{and}\quad
\mathbb{E}\xi_i\xi_j \le\mathbb{E}\zeta_i\zeta_j \qquad\mbox{for all $i\neq j$.}
\]
Let $F\dvtx\R^N \to\R$ be such that for all $\mathbf{x}\in\R^N$ and
$i\neq j$,
\[
\frac{\partial^2 F}{\partial x_i \,\partial x_j}(\mathbf{x}) \le0.
\]
Then $\mathbb{E} F(\bxi) \ge\mathbb{E} F(\bzeta)$.
\end{lemma}
\begin{pf*}{Proof of Theorem \ref{randomclassbound}}
Let $M\le N$ be fixed and choose $M$ sets from $\C$ uniformly at random
(without replacement). Let $\A$ denote the random subclass of
$\C$ obtained this way. Denote the likelihood ratio associated
to this class by
\[
L_\A(\bX) = \frac{{1}/{M} \sum_{S\in A} \phi_S(\bX)}{\phi
_0(\bX)}
= \frac{1}{M} \sum_{S\in\A} V_S,
\]
where $V_S= e^{\mu X_S - K\mu^2/2}$.
Then the optimal risk of the class $\C$ may be lower bounded by
\[
R^*_\C(\mu) - R^*_\A(\mu)
= \tfrac{1}{2} \bigl( \mathbb{E}_0| L_\A(\bX)-1 | - \mathbb{E}_0| L_\C(\bX)-1
| \bigr)
\ge- \tfrac{1}{2} \mathbb{E}_0 | L_\A(\bX) - L_\C(\bX) |.
\]
Denoting by $\widehat{\mathbb{E}}$ expectation with respect to the random
choice of $\A$, we have
\begin{eqnarray}
R^*_\C(\mu) & \ge& \widehat{\mathbb{E}} R^*_\A(\mu) -
\frac{1}{2} \mathbb{E}_0 \widehat{\mathbb{E}}
\biggl|\frac{1}{M} \sum_{S\in\A} V_S - \frac{1}{N} \sum_{S\in\C} V_S \biggr|
\nonumber\\
& \ge&
\widehat{\mathbb{E}} R^*_\A(\mu) -
\frac{1}{2} \sqrt{
\mathbb{E}_0 \widehat{\mathbb{E}}
\biggl(\frac{1}{M} \sum_{S\in\A} V_S - \frac{1}{N} \sum_{S\in\C} V_S \biggr)^2}
\nonumber\\
& \ge&
\widehat{\mathbb{E}} R^*_\A(\mu) -
\frac{1}{2} \sqrt{
\mathbb{E}_0 \biggl[\frac{1}{M} \cdot\frac{1}{N} \sum_{T\in\C}
\biggl( V_T - \frac{1}{N} \sum_{S\in\C} V_S \biggr)^2 \biggr] }
\nonumber\\
\eqntext{\mbox{(since the variance of a sample without replacement}} \\
\eqntext{\mbox{is less than that with replacement)}} \\
& = &
\widehat{\mathbb{E}} R^*_\A(\mu) -
\frac{1}{2\sqrt{M}}
\sqrt{
\frac{1}{N} \sum_{T\in\C}
\mathbb{E}_0 \biggl( V_T - \frac{1}{N} \sum_{S\in\C} V_S \biggr)^2}.\nonumber
\end{eqnarray}
An easy way to bound the
right-hand side is by writing
\begin{eqnarray*}
&&\mathbb{E}_0 \biggl( V_T - \frac{1}{N} \sum_{S\in\C} V_S \biggr)^2\\
&&\qquad \le
2\mathbb{E}_0 ( V_T - 1 )^2
+ 2\mathbb{E}_0 \biggl( 1 - \frac{1}{N} \sum_{S\in\C} V_S \biggr)^2 \\
&&\qquad \le
2\mathbb{E}_0 ( V_T - 1 )^2
+ \frac{2}{N} \sum_{S\in\C}\mathbb{E}_0 ( 1 - V_S )^2 \\
&&\qquad =
4 \Var(V_T) = 4 (e^{\mu^2K} -1 ).
\end{eqnarray*}
Summarizing, we have
\[
R^*_\C(\mu) \ge\widehat{\mathbb{E}} R^*_\A(\mu)-
\sqrt{\frac{e^{\mu^2K} -1}{M}} \ge\widehat{\mathbb{E}} R^*_\A(\mu)-
\frac{1}{4},
\]
where we used the assumption that $\mu\le\sqrt{(1/K)\log(M/16)}$.
Thus, it suffices to prove that
$\widehat{\mathbb{E}} R^*_\A(\mu) \ge1/2$.

We bound the optimal risk associated with $\A$ in terms of the
Bhattacharyya affinity
\[
\rho_\A(\mu) = \frac{1}{2}\mathbb{E}_0 \sqrt{\frac{(1/M)\sum_{S\in\A
} \phi
_S(\bX)}{\phi_0(\bX)}}
= \frac{1}{2} \mathbb{E}_0 \sqrt{\frac{1}{|\A|} \sum_{S\in\A} V_S}.
\]
Recalling from Section \ref{lowerbounds} that
$R^*_\A(\mu) \ge1-\sqrt{1-4 \rho_\A(\mu)^2}$ and using that
$\sqrt{1-4x^2}$ is concave, we have
\[
\widehat{\mathbb{E}} R^*_\A(\mu) \ge1- \sqrt{1-4 (\widehat{\mathbb
{E}}\rho_\A
(\mu) )^2}.
\]
Therefore, it suffices to show that the
expected Bhattacharyya affinity $\widehat{\mathbb{E}}\rho_\A(\mu) $
corresponding to the random class $\A$ satisfies
\[
\widehat{\mathbb{E}}\rho_\A(\mu)
= \frac{1}{2}\widehat{\mathbb{E}}\mathbb{E}_0 \sqrt{\frac{1}{|\A|} \sum
_{S\in
\A} V_S}
\ge\frac{\sqrt{3}}{4}.
\]
In the argument below, we fix the random class $\A$, relabel the
elements so that $\mathcal A=\{1,2,\ldots, |\mathcal A|\}$, and bound
$\rho_\A(\mu)$ from below.
Denote the minimum distance between any two elements of $\A$
by $\tau$.
To bound $\rho_\A(\mu)$,
we apply Slepian's lemma
with the function
\[
F(\mathbf{x}) = \sqrt{\frac1 {|\A|} \sum_{i=1}^{|\A|} e^{\mu
x_i-K\mu^2/2}},
\]
where $\mathbf{x}=(x_1, \ldots, x_{|\A|})$.
Simple calculation shows that the mixed second partial
derivatives of $F$ are negative, so Slepian's lemma is indeed applicable.

Next, we introduce the random vectors $\bxi$ and $\bzeta$.
Let the components of $\bxi$ be indexed by elements $S\in\A$ and define
$\xi_S= X_S=\sum_{i\in S} X_i$. Thus, under $\PROB_0$, each $\xi_S$ is
normal $(0,K)$
and $\mathbb{E} F(\bxi)$ is just the Bhattacharyya affinity $\rho_\A(\mu)$.
To define the random vector $\bzeta$, introduce $M+1$ independent
standard normal random variables: one variable $G_S$ for each $S \in\A$
and an extra variable $G_0$.
Recall that the definition of $\tau$ guarantees that the minimal
distance between any two elements of $\A$ as at least $\tau$.
Now let
\[
\zeta_S = G_S\frac{\tau}{\sqrt{2}} + G_0\sqrt{K-\frac{\tau^2}{2}}.
\]
Then clearly for each $S,T\in\A$, $\mathbb{E}\zeta_S^2 = K$ and
$\mathbb{E}\zeta_S\zeta_T = K-\tau^2/2$ ($S\neq T$).
On the other hand, $\mathbb{E}\xi_S^2= K$ and
\[
\mathbb{E}\xi_S\xi_T = |S\cap T| = K - \frac{d(S,T)^2}{2} \le K -\frac
{\tau^2}{2}
=\mathbb{E}\zeta_S\zeta_T.
\]
Therefore, by Slepian's lemma, $\rho_\A(\mu) = \mathbb{E} F(\bxi) \ge
\mathbb{E}
F(\bzeta)$.
However,
\begin{eqnarray*}
\mathbb{E} F(\bzeta) & = &
\mathbb{E}\sqrt{\frac{1}{|\A|} \sum_{S\in\A} e^{\mu\zeta_S-K\mu
^2/2}} \\
& = &
\mathbb{E}\sqrt{e^{\mu\sqrt{K-\tau^2/2}G_0-(K-\tau^2/2)\mu^2/2}
\frac{1}{|\A|} \sum_{S\in\A} e^{\mu\tau G_S/\sqrt{2}-\tau^2\mu^2/4}}
\\
& = &
\mathbb{E} e^{\mu\sqrt{K-\tau^2/2}G_0/2-(K-\tau^2/2)\mu^2/4}
\mathbb{E}\sqrt{ \frac{1}{|\A|} \sum_{S\in\A}
e^{\mu\tau G_S/\sqrt{2}-\tau^2\mu^2/4}} \\
& = &
e^{-\mu^2(K-\tau^2/2)/8}
\mathbb{E}\sqrt{ \frac{1}{|\A|}\sum_{S\in\A} e^{\mu\tau G_S/\sqrt
{2}-\tau
^2\mu^2/4}}.
\end{eqnarray*}
To finish the proof, it suffices to observe that the last expression
is the Bhattacharyya affinity corresponding to a class of disjoint sets,
all of size $\tau^2/2$, of cardinality $|\A|=M$. This case has been handled
in the first example of Section \ref{examples} where we showed that
\[
\mathbb{E}\sqrt{\frac{1}{|\A|} \sum_{S\in\A} e^{\mu\tau G_S/\sqrt
{2}-\tau
^2\mu^2/4}}
\ge R^*_\A\ge1 - \frac{1}{2}\sqrt{\frac{1}{M} e^{\mu^2 \tau^2/2}}
\ge\frac{3}{4},
\]
where again we used the condition $\mu\le\sqrt{\log(M/16)/K}$ and
the fact that $\tau^2/\break2 \le K$.

Therefore, under this condition on $\mu$, we have that for
any fixed $\A$,
\[
\rho_\A(\mu) = \tfrac{1}{2} \mathbb{E} F(\bzeta) \ge\tfrac{3}{8} e^{-\mu
^2(K-\tau^2/2)/8}
\]
and therefore
\[
\widehat{\mathbb{E}}\rho_\A(\mu) \ge\tfrac{3}{16}
e^{-\mu^2(K-t_\C(M)^2/2)/8},
\]
where $t_\C(M)$ is the median of $\tau$.
This concludes the proof.
\end{pf*}
\begin{remark*}[(An improvement)]
At the risk of losing a constant factor in the statement of
Theorem \ref{randomclassbound}, one may replace the parameter
$t_\C(M)$ by a larger quantity. The idea is that by thinning
the random subclass $\A$ one may consider a subset of $\A$
that has better separation properties. More precisely,
for an even integer $M\le N$ we may define a real-valued parameter
$\overline{t}_\C(M)>0$ of the class $\C$ as follows. Let $\A\subset
\C$
be obtained by choosing $M$ elements of $\C$ at random, without
replacement.
Order the elements $S_1,\ldots,S_M$ of $\A$ such that
\[
\min_{i\neq1} d(S_1,S_i) \ge\min_{i\neq2} d(S_2,S_i)
\ge\cdots\ge\min_{i\neq M} d(S_M,S_i)
\]
and define\vspace*{2pt} the subset $\widehat\A\subset\A$ by $\widehat\A=\{
A_1,\ldots
,A_{M/2}\}$.
Let the random variable $\overline\tau$ denote the smallest distance
between elements of $\widehat\A$ and let $\overline{t}_\C(M)$ be
the median
of $\tau$. It is easy to see that the proof of Theorem
\ref{randomclassbound} goes through, and one may replace
$t_\C(M)$ by $\overline{t}_\C(M)$ (by adjusting the constants appropriately).
One simply
needs to observe that
since each $V_S$ is
nonnegative,
\[
\rho_\A(\mu)
= \frac{1}{2} \mathbb{E}_0 \sqrt{\frac{1}{|\A|} \sum_{S\in\A} V_S}
\ge
\frac{1}{2} \mathbb{E}_0 \sqrt{\frac{1}{|\A|} \sum_{S\in\widehat\A} V_S}
= \frac{1}{\sqrt{2}} \rho_{\widehat\A}(\mu).
\]
If $\overline{t}_\C(M)$ is significantly larger than $t_\C(M)$, the
gain may be substantial.
\end{remark*}

If the class $\C$ is symmetric then thanks to Theorem
\ref{symmetricclass}, the theorem above can be improved and simplified.
If the class is symmetric, instead of having to work with
randomly chosen subclasses, one may optimally choose a separated
subset. Then the bounds can be
expressed in terms of the metric entropy of $\C$,
more precisely, by its
\textit{packing numbers} with respect to the canonical distance
$d(S,T) = \sqrt{\mathbb{E}_0 (X_S-X_T)^2}$.

We say that $\A\subset\C$ is a $t$-separated set (or $t$-packing)
if for any $S,T \in\A$, $d(S,T) \ge t$.
For $t <\sqrt{2K}$, define the \textit{packing number}
$M(t)$ as the size of a maximal
$t$-separated subset $\A$ of $\C$.
It is a simple well-known fact that packing numbers are closely
related to the covering numbers introduced in Section \ref{simpletests}
by the inequalities $N(t) \le M(t) \le N(t/2)$.
\begin{theorem}
Let $\C$ be symmetric in the sense of Theorem \ref{symmetricclass}
and let $t \le\sqrt{2K}$.
Then
\[
R_\C^*\ge1/2,
\]
whenever
\[
\mu\le\min\biggl(\sqrt{\frac{\log(M(t)/16)}{K}}, \frac{8\log(\sqrt
{3}/2)}{\sqrt{K-t^2/2}} \biggr).
\]
\end{theorem}
\begin{pf}
Let $\A\subset\C$ be a maximal $t$-separated subclass.
Since $\C$ is symmetric, by Theorem \ref{symmetricclass},
$R_\C^* \ge R_\A^*$ so it suffices to show that $R_\A^* \ge1/2$
for the indicated values of $\mu$. The rest of the proof is identical
to that of Theorem \ref{randomclassbound}.
\end{pf}

To interpret this result, take $t=\sqrt{2K(1-\varepsilon)}$ for
some $\varepsilon\in(0,1/2)$.
Then, by the theorem, $R^*\ge1/2$ if
\[
\mu\le\frac{1}{\sqrt{K}}
\min\biggl(\frac{8\log(\sqrt{3}/2)}{\sqrt{\varepsilon}},
\sqrt{\log\bigl(M\bigl(\sqrt{2K(1-\varepsilon)}\bigr)/16\bigr)}
\biggr).
\]
As an example, suppose that the class $\C$ is such that
there exists a constant $V>0$ such that
$M(t) \sim(n/t^2)^V$. (Recall that all classes with \textsc{vc} dimension
$V$ have an upper bound of this form for the packing numbers, see
remark on page \pageref{rmk:VC}.)
In this case, one may choose $\varepsilon^{-1} \sim V\log(n/K)$
and obtain that $R^*\ge1/2$ whenever
$\mu\le c \sqrt{(V/K)\log(n/K)}$ (for some constant $c$).
This closely matches the
bound obtained for the maximum test by Dudley's chaining bound.

\section{Optimal versus maximum test:
An analysis of the type I error}
\label{typeone}

In all examples considered above, upper bounds for the optimal risk
$R^*$ are derived by analyzing either the maximum test or the averaging
test. As the examples show, very often these simple tests have a
near-optimal performance. The optimal test $f^*$ is generally more
difficult to study. In this section, we analyze directly the
performance of
the optimal test. More precisely, we derive general upper bounds
for the type I error (i.e., the probability that the null hypothesis
is rejected under $\PROB_0$) of $f^*$. The upper bound involves
the expected value of the maximum of a Gaussian process indexed by a
sparse subset of $\C$ and can be significantly smaller than the
maximum over the whole class that appears in the performance bound
of the maximum test in Proposition \ref{maxtest}.
Unfortunately, we do not have an analogous bound for the type II error.

We consider the type I error of the optimal test $f^*$
\[
\PROB_0\{ f^*(\bX)=1 \}= \PROB_0\{ L(\bX)>1 \}
= \PROB_0 \biggl\{ \frac{1}{N} \sum_{S\in\C} e^{\mu X_S} > e^{K\mu^2/2}
\biggr\}.
\]
An easy bound is
$\frac{1}{N} \sum_{S\in\C} e^{\mu X_S} \le e^{\mu\max_{S\in\C}
X_S}$ so
\[
\PROB_0\{ L(\bX)>1 \} \le\PROB_0 \Bigl\{ \max_S X_S > K\mu/2 \Bigr\}.
\]
Thus, $\PROB_0\{ L(\bX)>1 \} \le\delta$ whenever
$\mu\ge(1/K)\mathbb{E}_0 \max_S X_S + \sqrt{(2/K)\log(1/\delta)}$.
Of course, we already know this from Proposition \ref{maxtest} where
this bound was derived for the (suboptimal) test based on maxima.

In order to understand the difference between the performance
of the optimal test $f^*$ and the maximum test,
one needs to compare the random variables
$(1/\mu)\log\frac{1}{N} \sum_{S\in\C} e^{\mu X_S}$ and
$\max_{S\in\C} X_S$.
\begin{proposition}
\label{typeonebound}
For any $\delta\in(0,1)$, the type \textup{I} error of the optimal
test $f^*$ satisfies
\[
\PROB_0\{ f^*(\bX)=1 \} \le\delta,
\]
whenever
\[
\mu\ge\frac{2}{K} \mathbb{E}_0 \max_{S \in\A} X_S+
\sqrt{\frac{32\log(2/\delta)}{K}},
\]
where $\A$ is any $\sqrt{K}/2$-cover of $\C$.
\end{proposition}

If $\A$ is a minimal $\sqrt{K}/2$-cover of $\C$, then
\[
(1/K) \mathbb{E}_0 \max_{S \in\A} X_S \le
\sqrt{\frac{2\log N(\sqrt{K}/2)}{K}}.
\]
By ``Sudakov's minoration'' [see
\citet{LeTa91}, Theorem 3.18] this upper bound is sharp up to a
constant factor.

It is instructive to compare this bound with that of
Proposition \ref{maxtest} for the performance of the maximum test.
In Proposition \ref{typeonebound}, we were able to replace
the expected maximum $\mathbb{E}_0 \max_{S\in\C} X_S$ by
$\mathbb{E}_0 \max_{S \in\A} X_S$ where now the maximum is taken over
a potentially much smaller subset $\A\subset\C$.
It is not difficult to construct examples when there is a substantial
difference, even in the order of magnitude, between the two
expected maxima so we have a genuine gain over the simple upper
bound of Proposition \ref{maxtest}.
Unfortunately, we do not
know if an analog upper bound holds for the type II error
$(1/N)\sum_{S\in\C} \PROB_S\{ f^*(\bX)=0 \}$ of the
optimal test $f^*$.
In cases when $\mathbb{E}_0 \max_{S \in\A} X_S \ll\mathbb{E}_0 \max_{S
\in
\C} X_S$,
we suspect that the maximum test is far from optimal.
However, to verify this conjecture, one would need a similar analysis
for the type II error as well.
\begin{pf*}{Proof of Proposition \ref{typeonebound}}
Introduce the notation
\[
M_\C(\mu)
= \mathbb{E}_0 \frac{1}{\mu} \log\biggl(\frac{1}{N} \sum_{S\in\C} e^{\mu
X_S} \biggr).
\]
Then
\begin{eqnarray*}
&&
\PROB_0 \biggl\{ \frac{1}{N} \sum_{S\in\C} e^{\mu X_S} > e^{K\mu^2/2}
\biggr\}
\\
&&\qquad =
\PROB_0 \biggl\{ \frac{1}{\mu}
\log\biggl(\frac{1}{N} \sum_{S\in\C} e^{\mu X_S} \biggr) > \frac{K\mu}{2} \biggr\}
\\
&&\qquad =
\PROB_0 \biggl\{ \frac{1}{\mu}
\log\biggl(\frac{1}{N} \sum_{S\in\C} e^{\mu X_S} \biggr)
-M_\C(\mu) > \frac{K\mu}{2} - M_\C(\mu) \biggr\}.
\end{eqnarray*}
We use Tsirelson's inequality (Lemma \ref{tsirelson}) to bound
this probability. To this end, we need to show that the
function $h\dvtx\R^N \to\R$ defined by
\[
h(\mathbf{x}) = \frac{1}{\mu}
\log\biggl(\frac{1}{N} \sum_{S\in\C} e^{\mu\sum_{i\in S} x_i} \biggr)
\]
is Lipschitz [where $\mathbf{x}=(x_1,\ldots,x_N)$].
Observing that
\[
\frac{\partial h}{\partial x_j} (\mathbf{x})
= \frac{{1/N} \sum_{S\in\C} \mathbh{1}_{\{ j\in S \}}e^{\mu
x_S} }
{{1/N} \sum_{S\in\C} e^{\mu x_S}} \in(0,1),
\]
we have
\[
\|\nabla h(\mathbf{x})\|^2 = \sum_{j=1}^n \biggl(\frac{\partial
h}{\partial x_j}(\mathbf{x}
) \biggr)^2
\le\sum_{j=1}^n \frac{\partial h}{\partial x_j} (\mathbf{x})
= K
\]
and therefore $h$ is indeed Lipschitz $\sqrt{K}$. By Tsirelson's
inequality, we have
\[
\PROB_0\{ f^*(\bX)=1 \} \le
\exp\biggl( - \frac{ (K\mu/2 - M_\C(\mu) )^2}{2K} \biggr).
\]
Thus, the type I error is bounded by $\delta$ if
\[
\mu\ge\frac{2M_{\C}(\mu)}{K} + \sqrt{\frac{8}{K}\log\frac
{1}{\delta}}.
\]
It remains to bound $M_{\C}(\mu)$.

Let $t\le\sqrt{2K}$ be a positive integer and
consider a minimal $t$-cover of the set~$\C$, that is, a set
$\A\subset\C$ with cardinality $|\A| = N(t)$ such that, if
$\pi(S)$ denotes an element in $\A$ whose distance to $S\in\C$
is minimal then $d(S,\pi(S)) \le t$ for all $S\in\C$.
Then clearly,
\[
M_{\C}(\mu) \le
\mathbb{E}_0 \frac{1}{\mu}
\log\biggl(\frac{1}{N} \sum_{S\in\C} e^{\mu(X_S-X_{\pi(S)})} \biggr)
+ \mathbb{E}_0 \max_{S \in\A} X_S.
\]
To bound the first term on the right-hand side, note that, by
Jensen's inequality,
\[
\mathbb{E}_0 \frac{1}{\mu}
\log\biggl(\frac{1}{N} \sum_{S\in\C} e^{\mu(X_S-X_{\pi(S)})} \biggr)
\le\frac{1}{\mu}
\log\biggl(\frac{1}{N} \sum_{S\in\C} \mathbb{E}_0 e^{\mu(X_S-X_{\pi(S)})} \biggr)
\le\frac{\mu t^2 }{2}
\]
since for each $S$, $d_H(X_S,X_{\pi(S)}) \le t^2$ and therefore
$X_S-X_{\pi(S)}$ is a centered normal random variable with variance
$d_H(X_S,X_{\pi(S)})$.
For the second term, we have
\[
\mathbb{E}_0 \max_{S \in\A} X_S \le\sqrt{2K\log N(t)}.
\]
Choosing $t^2=K/4$, we obtain the proposition.
\end{pf*}

\section*{Acknowledgments}
We thank Ery Arias-Castro and Emmanuel Cand\`es
for discussions on the topic of the paper.
We also thank the referees for their valuable remarks.

\printaddresses

\end{document}